\documentclass{article}%
\usepackage{amsfonts}
\usepackage{graphicx}
\usepackage{amsmath}
\usepackage{amssymb}%
\newtheorem{theorem}{Theorem}[section]
\newtheorem{definition}[theorem]{Definition}
\newtheorem{lemma}[theorem]{Lemma}
\newtheorem{remark}[theorem]{Remark}
\newenvironment{proof}[1][Proof]{\textbf{#1.} }{\ \rule{0.5em}{0.5em}}
\begin{document}

\title{Transfinite thin plate spline interpolation}
\date{8 October 2009}
\author{Aurelian Bejancu\\Department of Mathematics, Kuwait University\\PO Box 5969, Safat 13060, Kuwait }
\maketitle

\begin{abstract}
Duchon's method of thin plate splines defines a polyharmonic interpolant to
scattered data values as the minimizer of a certain integral functional. For
transfinite interpolation, i.e.\ interpolation of continuous data prescribed
on curves or hypersurfaces, Kounchev has developed the method of polysplines,
which are piecewise polyharmonic functions of fixed smoothness across the
given hypersurfaces and satisfy some boundary conditions. Recently, Bejancu
has introduced boundary conditions of Beppo Levi type to construct a
semi-cardinal model for polyspline interpolation to data on an infinite set of
parallel hyperplanes. The present paper proves that, for periodic data on a
finite set of parallel hyperplanes, the polyspline interpolant satisfying
Beppo Levi boundary conditions is in fact a thin plate spline, i.e.\ it
minimizes a Duchon type functional.\medskip

\noindent\emph{MSC:} 41A05; 41A15; 41A63 \medskip

\noindent\emph{Keywords:} Transfinite interpolation; Polyharmonic functions;
Boundary conditions; Variational principle; Polysplines; Natural $\mathcal{L}%
$-splines; Radial basis functions

\end{abstract}


\section{Introduction}

\label{intro}

Thin plate spline interpolation, one of the main algorithms for multivariable
scattered data approximation, was introduced by Duchon in \cite{Du76}%
-\cite{Du77} and has since generated a vast amount of research on its theory
and applications (for a comprehensive survey, see Wendland \cite{We05}). Among
all functions $F:\mathbb{R}^{d}\rightarrow\mathbb{R}$ having all partial
derivatives $\partial^{\alpha}F$ of total order $\left\vert \alpha\right\vert
=p$ in $L^{2}\left(  \mathbb{R}^{d}\right)  $ (where $p>d/2$) and taking
prescribed values at a finite (sufficiently large) number of fixed scattered
points in $\mathbb{R}^{d}$, the \emph{thin plate spline} (or \emph{surface
spline}) interpolant to the data is defined as the unique minimizer of the
integral functional%
\begin{equation}
\int\limits_{\mathbb{R}^{d}}\sum_{\left\vert \alpha\right\vert =p}\frac
{p!}{\alpha!}\left\vert \partial^{\alpha}F\left(  x\right)  \right\vert
^{2}dx.\label{eq:Du}%
\end{equation}
Duchon's theory identifies this minimizer as a function $S\in C^{2p-d-1}%
\left(  \mathbb{R}^{d}\right)  $ which is polyharmonic of order $p$, i.e.
$\Delta^{p}S\left(  x\right)  =0$ with $\Delta$ the Laplace operator, for all
$x\in\mathbb{R}^{d}$ except the given scattered points. The computational
usefulness of $S$ stems from its explicit radial basis representation as a
finite linear combination of translates of the fundamental solution of
$\Delta^{p}$ in $\mathbb{R}^{d}$.

A different interpolation problem that can be formulated in the multivariable
case is that of constructing a surface or function that matches continuous
data prescribed on some collection of curves or hypersurfaces. In computer
aided design, this problem is referred to as \emph{transfinite interpolation}
(Sabin \cite{MAS96}). Typical examples include the reconstruction of 2D
surfaces from level curves or from track data, as well as the visualization of
3D objects from scan data. A variational solution for transfinite
interpolation along a set of curves in $\mathbb{R}^{2}$ is defined by Apprato
and Arcang\'{e}li \cite{AA91} (see also \cite[Chapter~X]{ALT04}), and its
approximate computation is obtained by smoothing finite elements. On the other
hand, Kounchev's polyspline method treated in the monograph \cite{Ku01} is
based on explicit piecewise polyharmonic functions of any number of variables.

The present paper studies the polyspline method in a simplified setting,
namely for periodic data given on parallel hyperplanes. Let $n$, $p$, $N$ be
positive integers with $2\leq p\leq N+1$, and $\tau:=\left\{  t_{0}%
,\ldots,t_{N}\right\}  $, where $t_{0}<t_{1}<\ldots<t_{N}$ are fixed real
numbers. Let $\Omega$ be the closure in $\mathbb{R}^{n+1}$ of the union of the
open strips%
\begin{equation}
\Omega_{j}:=\left\{  \left(  t,y\right)  \in\mathbb{R}\times\mathbb{R}%
^{n}:t\in\left(  t_{j-1},t_{j}\right)  \right\}  ,\quad j\in\left\{
1,2,\ldots,N\right\}  .\label{eq:strips}%
\end{equation}
A function $S:\Omega\rightarrow\mathbb{C}$ is called a \emph{polyspline} of
order $p$ on strips determined by $\tau$ in $\Omega$ if $S\in C^{2p-2}\left(
\Omega\right)  $, $S$ is polyharmonic of order $p$ on each open strip
$\Omega_{j}$, $j\in\left\{  1,\ldots,N\right\}  $, and $S$ is $2\pi$-periodic
in each of its last $n$ variables. In order to determine a unique polyspline
that interpolates continuous periodic data prescribed on the hyperplanes
$\left\{  t_{j}\right\}  \times\mathbb{R}^{n}$, $j\in\left\{  0,1,\ldots
,N\right\}  $, additional boundary conditions are imposed on $\left\{
t_{0}\right\}  \times\mathbb{R}^{n}$ and $\left\{  t_{N}\right\}
\times\mathbb{R}^{n}$. Then, provided that $p$ is an even integer, Kounchev's
theory \cite[Theorem~20.14]{Ku01} identifies the polyspline interpolant as
minimizer of the functional%
\begin{equation}
\int\limits_{\left[  t_{0},t_{N}\right]  \times\mathbb{T}^{n}}\left\vert
\Delta^{p/2}F\left(  t,y\right)  \right\vert ^{2}dtdy,\label{eq:Ku}%
\end{equation}
subject to the boundary and interpolation conditions, where $\Delta$ is now
the Laplace operator in $\mathbb{R}^{n+1}$, and $\mathbb{T}:=\left[  -\pi
,\pi\right]  $.

Recently, Bejancu \cite{AB08} has introduced boundary conditions of Beppo Levi
type in the construction of semi-cardinal polyspline interpolation to data on
the infinite set of parallel hyperplanes $\left\{  j\right\}  \times
\mathbb{R}^{n}$, $j\in\left\{  0,1,\ldots\right\}  $. To incorporate such
conditions in the above setting, let $t_{-1}:=-\infty$, $t_{N+1}:=\infty$, and
let the corresponding open strips $\Omega_{0}$, $\Omega_{N+1}$ be defined as
in (\ref{eq:strips}).

\begin{definition}
\label{def:poly-nat} Let $p\geq2$ be a fixed integer. A function
$S:\mathbb{R}^{n+1}\rightarrow\mathbb{C}$ is called a \emph{Beppo Levi
polyspline of order }$p$\emph{\ on strips determined by }$\tau$\emph{\ in
}$R^{n+1}$ if the following conditions hold:

\noindent(i) $S\in C^{2p-2}\left(  \mathbb{R}^{n+1}\right)  $;

\noindent(ii) $S$ is polyharmonic of order $p$ on each open strip $\Omega_{j}%
$, $j\in\left\{  0,1,\ldots,N+1\right\}  $;

\noindent(iii) $S$ is $2\pi$-periodic in each of its last $n$ variables;

\noindent(iv) $S$ satisfies the Beppo Levi conditions%
\begin{equation}
\partial^{\alpha}S\in L^{2}\left(  \mathbb{R}\times\mathbb{T}^{n}\right)
,\quad\forall\left\vert \alpha\right\vert =p.\label{eq:beppolevi}%
\end{equation}

\noindent The space of all Beppo Levi polysplines of order $p$ on strips
determined by $\tau$ in $\mathbb{R}^{n+1}$ is denoted by $\mathcal{S}%
_{p}\left(  \tau,n\right)  $.
\end{definition}

\noindent The main results of the present paper (obtained in section~3) prove
that, for periodic data functions prescribed on $\left\{  t_{j}\right\}
\times\mathbb{R}^{n}$, $j\in\left\{  0,1,\ldots,N\right\}  $, there exists a
unique interpolant in $\mathcal{S}_{p}\left(  \tau,n\right)  $, which
minimizes the Duchon type functional%
\begin{equation}
\int\limits_{\mathbb{R}\times\mathbb{T}^{n}}\sum_{\left\vert \alpha\right\vert
=p}\frac{p!}{\alpha!}\left\vert \partial^{\alpha}F\left(  t,y\right)
\right\vert ^{2}dtdy.\label{eq:Du-per}%
\end{equation}
This shows that the Beppo Levi polyspline interpolant is a genuine thin plate
spline analog for transfinite interpolation. Moreover, our variational
characterization is valid for any order $p\geq2$, without the restriction that
$p>d/2$ as in (\ref{eq:Du}) or that $p$ is even as in (\ref{eq:Ku}). In the
excluded case $p=1$, the results degenerate to the classical existence,
uniqueness, and variational characterization of harmonic solutions to the
Dirichlet problem in each strip separately.

Note that the periodicity assumption is of importance for both theoretical and
practical purposes. In our proofs, it will serve to simplify certain technical
arguments that guarantee the uniqueness of the polyspline constructions. On
the other hand, periodic polysplines are directly applicable to the problem of
transfinite interpolation on a cylinder domain, as demonstrated, for example,
in the case of visualization of the heart surface in medical imaging
\cite{KW03}.

The construction of periodic polysplines on strips is reduced to that of a
family of exponential $\mathcal{L}$-splines by separation of variables. To
describe this procedure in our setting, let $S\in\mathcal{S}_{p}\left(
\tau,n\right)  $ and, for each $t\in\mathbb{R}$, define the sequence of
Fourier coefficients of $S$ with respect to its last $n$ variables by%
\[
\widehat{S}_{\xi}\left(  t\right)  :=\int_{\mathbb{T}^{n}}e^{-i\left\langle
\xi,y\right\rangle }S\left(  t,y\right)  dy,\quad\xi\in\mathbb{Z}^{n},
\]
where $\left\langle \xi,y\right\rangle $ is the dot product in $\mathbb{R}%
^{n}$. For a fixed $\xi\in\mathbb{Z}^{n}$, note that condition \textit{(i)}
above implies $\widehat{S}_{\xi}\in C^{2p-2}\left(  \mathbb{R}\right)  $.
Further, let $\left\vert \xi\right\vert $ the Euclidean norm of $\xi$ in
$\mathbb{R}^{n}$ and consider the ordinary differential operator%
\begin{equation}
\mathcal{L}_{\xi}:=\left(  \frac{d^{2}}{dt^{2}}-\left\vert \xi\right\vert
^{2}\right)  ^{p},\label{eq:Lxi}%
\end{equation}
with null-space%
\[
\mathrm{Ker}\mathcal{L}_{\xi}=\left\{
\begin{array}
[c]{ll}%
\mathrm{span}\left\{  e^{\pm\left\vert \xi\right\vert t},te^{\pm\left\vert
\xi\right\vert t},\ldots,t^{p-1}e^{\pm\left\vert \xi\right\vert t}\right\}
, & \text{if }\xi\not =0\text{,}\\
\mathrm{span}\left\{  1,t,\ldots,t^{2p-1}\right\}  , & \text{if }\xi=0\text{.}%
\end{array}
\right.
\]
From the polyharmonic condition \textit{(ii)}, we deduce (as in
\cite[Lemma~2.1]{BKR07})
\[
\mathcal{L}_{\xi}\widehat{S}_{\xi}\left(  t\right)  =0,\quad\forall
t\in\left(  t_{j-1},t_{j}\right)  ,\;\forall j\in\left\{  0,1,\ldots
,N+1\right\}  .
\]
Also, the Beppo Levi conditions (\ref{eq:beppolevi}) imply (as in \cite[Eq.
(4)]{AB08})%
\[
\widehat{S}_{\xi}\in\left\{
\begin{array}
[c]{l}%
\mathrm{Ker}\left(  \frac{d}{dt}-\left\vert \xi\right\vert \right)  ^{p}%
\quad\text{on the interval }\left(  -\infty,t_{0}\right)  ,\\
\mathrm{Ker}\left(  \frac{d}{dt}+\left\vert \xi\right\vert \right)  ^{p}%
\quad\text{on the interval }\left(  t_{N},\infty\right)  .
\end{array}
\right.
\]
These observations motivate the following definition.

\begin{definition}
\label{def:natLspline}For a fixed $\xi\in\mathbb{Z}^{n}$, the function
$s:\mathbb{R}\rightarrow\mathbb{C}$ is called a \emph{natural }$\mathcal{L}%
_{\xi}$\emph{-spline on }$\tau$ if:

\noindent(i) $s\in C^{2p-2}\left(  \mathbb{R}\right)  $;

\noindent(ii) $\mathcal{L}_{\xi}s\left(  t\right)  =0,\quad\forall t\in\left(
t_{j-1},t_{j}\right)  ,\;\forall j\in\left\{  0,1,\ldots,N+1\right\}  $;

\noindent(iii) $\left(  \frac{d}{dt}-\left\vert \xi\right\vert \right)
^{p}s\left(  t\right)  =0$, $\forall t<t_{0}$, and $\left(  \frac{d}%
{dt}+\left\vert \xi\right\vert \right)  ^{p}s\left(  t\right)  =0$, $\forall
t>t_{N}$.

\noindent The space of all natural $\mathcal{L}_{\xi}$-splines on $\tau$ will
be denoted by $\mathcal{S}_{p,\xi}\left(  \tau\right)  $.
\end{definition}

\begin{remark}
\emph{ For $\xi=0$, this definition corresponds to the well-known natural
polynomial splines of degree $2p-1$. For $\xi\not =0$, the use of adjoint
operators on the two extreme intervals in condition (iii) is equivalent to the
decay of the natural $\mathcal{L}_{\xi}$-spline $s$ at $\pm\infty$. This
differs from the standard natural conditions for Chebyshev splines (cf.
Schumaker \cite[p.\ 396]{LLS81}), which employ one and the same natural
operator on both sides. It is of interest to note that in the case of
semi-cardinal interpolation with $\mathcal{L}_{\xi}$-splines treated in
\cite{AB08}, the splitting of the operator $\mathcal{L}_{\xi}$ into its two
adjoint factors on the left/right boundary pieces is a direct consequence of
the Wiener-Hopf factorization technique. }
\end{remark}

The above arguments show that $S\in\mathcal{S}_{p}\left(  \tau,n\right)  $
implies $\widehat{S}_{\xi}\in\mathcal{S}_{p,\xi}\left(  \tau\right)  $,
$\forall\xi\in\mathbb{Z}^{n}$. Conversely, our results on Beppo Levi
polyspline interpolation in section~3 will follow from the properties of the
natural exponential $\mathcal{L}_{\xi}$-splines. The construction and
necessary analysis of the natural $\mathcal{L}_{\xi}$-spline interpolation
schemes are contained in section~2. An original contribution of this analysis
is the use of radial basis representations to estimate the effect of the
variable parameter $\xi$ on the size of the Lagrange functions for natural
$\mathcal{L}_{\xi}$-spline interpolation.

The extension of our results to the setting of data prescribed on concentric
spheres will have to address significant technical differences. A further
problem of interest would be to study the convergence properties of Beppo Levi
polyspline interpolation.

\section{Natural $\mathcal{L}_{\xi}$-spline interpolation}

\label{L-splines} Note that the operator $\mathcal{L}_{\xi}$ of (\ref{eq:Lxi})
and the space $\mathcal{S}_{p,\xi}\left(  \tau\right)  $ of natural
$\mathcal{L}_{\xi}$-splines can be defined not only for $\xi\in\mathbb{Z}^{n}%
$, but for any $\xi\in\mathbb{R}^{n}$. Therefore in this section we work with
$\xi\in\mathbb{R}^{n}\backslash\left\{  0\right\}  $, the version of the
results for $\xi=0$ being well-known. We will also assume that $p\geq2$
throughout the section.

Our first result asserts existence and uniqueness for the problem of natural
$\mathcal{L}_{\xi}$-spline interpolation at the set of knots $\tau=\left\{
t_{0},t_{1},\ldots,t_{N}\right\}  $.

\begin{theorem}
\label{thm:EU-natLspline}For any set of data values $\left\{  y_{0}%
,\ldots,y_{N}\right\}  \subset\mathbb{R}$, there exists a unique natural
$\mathcal{L}_{\xi}$-spline on $\tau$, $s\in\mathcal{S}_{p,\xi}\left(
\tau\right)  $, such that%
\begin{equation}
s\left(  t_{j}\right)  =y_{j},\quad j\in\left\{  0,1,\ldots,N\right\}
.\label{eq:1D-intconds}%
\end{equation}

\end{theorem}

\noindent We will derive the proof of this theorem from the following result.

\begin{theorem}
\label{thm:FI-natLspline}Let $\sigma\in\mathcal{S}_{p,\xi}\left(  \tau\right)
$ be any natural $\mathcal{L}_{\xi}$-spline on $\tau$ and let $\psi\in
C^{p}\left(  \mathbb{R}\right)  $ such that $\psi^{\left(  m\right)  }\in
L^{2}\left(  \mathbb{R}\right)  $, $\forall m\in\left\{  0,\ldots,p\right\}
$, and $\psi\left(  t_{j}\right)  =0$, $\forall j\in\left\{  0,\ldots
,N\right\}  $. Then the following `fundamental identity'\ holds:%
\begin{equation}
\int_{-\infty}^{\infty}\left(  \frac{d}{dt}-\left\vert \xi\right\vert \right)
^{p}\sigma\left(  t\right)  \left(  \frac{d}{dt}-\left\vert \xi\right\vert
\right)  ^{p}\overline{\psi\left(  t\right)  }dt=0.\label{eq:FI-natLspline}%
\end{equation}

\end{theorem}

\noindent Another direct consequence of (\ref{eq:FI-natLspline}) is the
variational characterization of the natural $\mathcal{L}_{\xi}$-spline $s$ of
Theorem~\ref{thm:EU-natLspline}. Although this will not be employed further,
it may be of interest to mention it here.

\begin{theorem}
\label{thm:var-p1D}Given the set of values $\left\{  y_{0},\ldots
,y_{N}\right\}  \subset\mathbb{R}$, let $s$ be the unique natural
$\mathcal{L}_{\xi}$-spline on $\tau$ satisfying the interpolation conditions
(\ref{eq:1D-intconds}). If $f\in C^{p}\left(  \mathbb{R}\right)  $ is any
other function such that $f^{\left(  m\right)  }\in L^{2}\left(
\mathbb{R}\right)  $, $\forall m\in\left\{  0,\ldots,p\right\}  $, and
$f\left(  t_{j}\right)  =y_{j} $, $\forall j\in\left\{  0,\ldots,N\right\}  $,
then%
\[
\int_{-\infty}^{\infty}\left\vert \left(  \frac{d}{dt}-\left\vert
\xi\right\vert \right)  ^{p}s\left(  t\right)  \right\vert ^{2}dt<\int
_{-\infty}^{\infty}\left\vert \left(  \frac{d}{dt}-\left\vert \xi\right\vert
\right)  ^{p}f\left(  t\right)  \right\vert ^{2}dt.
\]

\end{theorem}

A related question at this point is whether the statements of the last two
theorems remain true if the left natural operator $\left(  \frac{d}%
{dt}-\left\vert \xi\right\vert \right)  ^{p}$ is replaced by its adjoint
$\left(  -\frac{d}{dt}-\left\vert \xi\right\vert \right)  ^{p}$. One way to
see that the answer is positive is to go through the proof of
Theorem~\ref{thm:FI-natLspline} in subsection~2.1 and make the requisite
changes. A faster way, however, is to use the next result that will also play
a significant role in the arguments of section~3.

\begin{lemma}
\label{le:identity}Let $\sigma\in\mathcal{S}_{p,\xi}\left(  \tau\right)  $ be
any natural $\mathcal{L}_{\xi}$-spline on $\tau$ and let $\psi\in C^{p}\left(
\mathbb{R}\right)  $ such that $\psi^{\left(  m\right)  }\in L^{2}\left(
\mathbb{R}\right)  $, $\forall m\in\left\{  0,\ldots,p\right\}  $. Then%
\begin{align}
& \int_{-\infty}^{\infty}\left(  \frac{d}{dt}-\left\vert \xi\right\vert
\right)  ^{p}\sigma\left(  t\right)  \left(  \frac{d}{dt}-\left\vert
\xi\right\vert \right)  ^{p}\overline{\psi\left(  t\right)  }%
dt\label{eq:identity}\\
& =\sum_{m=0}^{p}\left(
\begin{array}
[c]{c}%
p\\
m
\end{array}
\right)  \left\vert \xi\right\vert ^{2\left(  p-m\right)  }\int_{-\infty
}^{\infty}\sigma^{\left(  m\right)  }\left(  t\right)  \,\overline
{\psi^{\left(  m\right)  }\left(  t\right)  }\,dt,\nonumber
\end{align}
where $\left(
\begin{array}
[c]{c}%
p\\
m
\end{array}
\right)  =\frac{p!}{m!\left(  p-m\right)  !}$ is the usual binomial coefficient.
\end{lemma}

\noindent The proof of this identity in subsection~2.2 works by showing that,
after expanding the two brackets in the left-hand side, integrating term by
term, and collecting terms of the same power of $\left\vert \xi\right\vert $,
the odd power terms vanish. Since replacing $\frac{d}{dt}-\left\vert
\xi\right\vert $ by $\frac{d}{dt}+\left\vert \xi\right\vert $ in the left-hand
side will only change the sign of the odd power terms, we deduce that, under
the hypotheses of Lemma~\ref{le:identity},%
\begin{align*}
& \int_{-\infty}^{\infty}\left(  \frac{d}{dt}+\left\vert \xi\right\vert
\right)  ^{p}\sigma\left(  t\right)  \left(  \frac{d}{dt}+\left\vert
\xi\right\vert \right)  ^{p}\overline{\psi\left(  t\right)  }\, dt\\
& =\int_{-\infty}^{\infty}\left(  \frac{d}{dt}-\left\vert \xi\right\vert
\right)  ^{p}\sigma\left(  t\right)  \left(  \frac{d}{dt}-\left\vert
\xi\right\vert \right)  ^{p}\overline{\psi\left(  t\right)  }\, dt.
\end{align*}
Therefore Theorems \ref{thm:FI-natLspline} and \ref{thm:var-p1D} also hold
with $\frac{d}{dt}+\left\vert \xi\right\vert $ in place of $\frac{d}%
{dt}-\left\vert \xi\right\vert $.

It can be noticed that, for a fixed $\xi$, Theorems~\ref{thm:EU-natLspline},
\ref{thm:FI-natLspline}, and \ref{thm:var-p1D} follow along a rather classical
route in spline theory, the only novelty being due to the adjoint boundary
operators in Definition~\ref{def:natLspline} of the natural $\mathcal{L}_{\xi
}$-spline. The classical theory, however, does not cover the dependence of the
natural $\mathcal{L}_{\xi}$-spline interpolation scheme on the variable
parameter $\xi$, as will be needed in section~3. In order to study this
dependence, for each $\xi\in\mathbb{Z}^{n}$ and $j\in\left\{  0,1,\ldots
,N\right\}  $, let $L_{\xi,j}\in\mathcal{S}_{p,\xi}\left(  \tau\right)  $ be
the unique \emph{Lagrange function} determined by the interpolation conditions%
\begin{equation}
L_{\xi,j}\left(  t_{j}\right)  =1\text{ and }L_{\xi,j}\left(  t_{k}\right)
=0\text{ for }k\in\left\{  0,1,\ldots,N\right\}  \backslash\left\{  j\right\}
.\label{eq:Lagrange-int}%
\end{equation}
Hence, if $s$ is the natural $\mathcal{L}_{\xi}$-spline that interpolates the
values $y_{0}$, $y_{1}$, $\ldots$ , $y_{N}$ in Theorem~\ref{thm:EU-natLspline}%
, we have the Lagrange formula%
\begin{equation}
s\left(  t\right)  =\sum_{j=0}^{N}y_{j}L_{\xi,j}\left(  t\right)
,\quad\forall t\in\mathbb{R}.\label{eq:Lagrange}%
\end{equation}
The following theorem estimates the effect of the parameter $\xi$ on the
stability of the Lagrange scheme (\ref{eq:Lagrange}).

\begin{theorem}
\label{thm:stability}There exists a constant $C_{0}=C_{0}\left(
p,\tau\right)  >0$ such that, for all $j\in\left\{  0,\ldots,N\right\}  $,
$m\in\left\{  0,..,2p-2\right\}  $, and all $\xi\in\mathbb{R}^{n}$ with
$\left\vert \xi\right\vert \geq\frac{1}{2}$, we have:
\begin{equation}
\left\vert \frac{d^{m}}{dt^{m}}L_{\xi,j}\left(  t\right)  \right\vert \leq
C_{0}\left(  1+\left\vert \xi\right\vert ^{m}\right)  ,\quad\forall
t\in\mathbb{R}.\label{eq:stability}%
\end{equation}

\end{theorem}

A similar statement for $m=0$ and $t\in\left[  t_{0},t_{N}\right]  $ was given
by Kounchev \cite[Lemma~3]{Ku94} in the context of an $\mathcal{L}_{\xi}%
$-spline interpolation scheme with different boundary conditions and an even
integer $p$. The validity of the present result for all orders $m\in\left\{
0,..,2p-2\right\}  $ is of crucial importance for the analysis of section~3.
In subsection~2.3 we will employ an original method of proof of this result
via `radial basis' representations of the Lagrange functions $L_{\xi,j}$.

\subsection{Proofs of Theorems~\ref{thm:EU-natLspline},
\ref{thm:FI-natLspline}, and \ref{thm:var-p1D}}

We start by establishing the following auxiliary result.

\begin{lemma}
\label{le:pol-growth}If $\psi\in C^{1}\left(  \mathbb{R}\right)  $ and
$\psi^{\prime}\in L^{2}\left(  \mathbb{R}\right)  $, then there exists
$C_{\psi}\geq0$ such that%
\[
\left\vert \psi\left(  t\right)  \right\vert \leq C_{\psi}\left(  1+\left\vert
t\right\vert ^{1/2}\right)  ,\quad\forall t\in\mathbb{R}.
\]

\end{lemma}

\noindent\proof Combining the Leibniz-Newton formula $\psi\left(  t\right)
=\psi\left(  0\right)  +\int_{0}^{t}\psi^{\prime}\left(  u\right)  du$ with
the Cauchy-Schwarz inequality, we obtain for $t>0$:%
\begin{align*}
\left\vert \psi\left(  t\right)  \right\vert  & \leq\left\vert \psi\left(
0\right)  \right\vert +\left(  \int_{0}^{t}du\right)  ^{1/2}\left(  \int
_{0}^{t}\left\vert \psi^{\prime}\left(  u\right)  \right\vert ^{2}du\right)
^{1/2}\\
& \leq\left\vert \psi\left(  0\right)  \right\vert +t^{1/2}\left(
\int_{\mathbb{R}}\left\vert \psi^{\prime}\left(  u\right)  \right\vert
^{2}du\right)  ^{1/2}.
\end{align*}
The conclusion follows by making a similar estimate for $t<0$ and letting
$C_{\psi}:=\max\left\{  \left\vert \psi\left(  0\right)  \right\vert
,\left\Vert \psi^{\prime}\right\Vert _{L^{2}\left(  \mathbb{R}\right)
}\right\}  $. \endproof \medskip

\noindent\textbf{Proof of Theorem \ref{thm:FI-natLspline}.} Let $J_{\xi}$
denote the convergent integral of (\ref{eq:FI-natLspline}). Since $\left(
\frac{d}{dt}-\left\vert \xi\right\vert \right)  ^{p}\sigma\left(  t\right)  =0
$ for $t\leq t_{0}$, the integration domain of $J_{\xi}$ can be replaced by
the interval $[t_{0},\infty)$. Using integration by parts,%
\begin{align*}
J_{\xi}  & =\left[  \left(  \frac{d}{dt}-\left\vert \xi\right\vert \right)
^{p}\sigma\left(  t\right)  \left(  \frac{d}{dt}-\left\vert \xi\right\vert
\right)  ^{p-1}\overline{\psi\left(  t\right)  }\right]  _{t=t_{0}}^{\infty}\\
& -\int_{t_{0}}^{\infty}\left(  \frac{d}{dt}+\left\vert \xi\right\vert
\right)  \left(  \frac{d}{dt}-\left\vert \xi\right\vert \right)  ^{p}%
\sigma\left(  t\right)  \left(  \frac{d}{dt}-\left\vert \xi\right\vert
\right)  ^{p-1}\overline{\psi\left(  t\right)  }dt.
\end{align*}
The expression in square brackets vanishes at $t_{0}$ due to the above
condition on $\sigma$, and it also vanishes at $\infty$ due to the exponential
decay of $\sigma$ at $\infty$ and the fact that, by Lemma~\ref{le:pol-growth},
the function $\left(  \frac{d}{dt}-\left\vert \xi\right\vert \right)
^{p-1}\overline{\psi\left(  t\right)  }$ has at most an algebraic growth at
$\infty$. Applying similar arguments successively,%
\[
J_{\xi}=\left(  -1\right)  ^{p-1}\int_{t_{0}}^{\infty}\left(  \frac{d}%
{dt}+\left\vert \xi\right\vert \right)  ^{p-1}\left(  \frac{d}{dt}-\left\vert
\xi\right\vert \right)  ^{p}\sigma\left(  t\right)  \left(  \frac{d}%
{dt}-\left\vert \xi\right\vert \right)  \overline{\psi\left(  t\right)  }dt.
\]
Since $\sigma\in\mathcal{S}_{p,\xi}\left(  \tau\right)  $, for each
$j\in\left\{  0,\ldots,N\right\}  $ there exists a constant $\sigma_{j}$ such
that%
\begin{align*}
&  \left(  \frac{d}{dt}+\left\vert \xi\right\vert \right)  ^{p-1}\left(
\frac{d}{dt}-\left\vert \xi\right\vert \right)  ^{p}\sigma\left(  t\right) \\
&  =\sigma_{j}\left(  \frac{d}{dt}+\left\vert \xi\right\vert \right)
^{p-1}\left(  \frac{d}{dt}-\left\vert \xi\right\vert \right)  ^{p}\left(
t^{p-1}e^{-\left\vert \xi\right\vert t}\right) \\
&  =\sigma_{j}\left(  p-1\right)  !\left(  -2\left\vert \xi\right\vert
\right)  ^{p}e^{-\left\vert \xi\right\vert t},\quad\forall t\in\left(
t_{j},t_{j+1}\right)  ,
\end{align*}
where $t_{N+1}:=\infty$. Hence%
\[
J_{\xi}=-\left(  p-1\right)  !\left(  2\left\vert \xi\right\vert \right)
^{p}\sum_{j=0}^{N}\sigma_{j}\int_{t_{j}}^{t_{j+1}}e^{-\left\vert
\xi\right\vert t}\left(  \frac{d}{dt}-\left\vert \xi\right\vert \right)
\overline{\psi\left(  t\right)  }dt.
\]
Noting that $e^{-\left\vert \xi\right\vert t}\left(  \frac{d}{dt}-\left\vert
\xi\right\vert \right)  \overline{\psi\left(  t\right)  }=\frac{d}{dt}\left(
e^{-\left\vert \xi\right\vert t}\overline{\psi\left(  t\right)  }\right)  $
and using the hypotheses $\psi\left(  t_{j}\right)  =0$, $j\in\left\{
0,1,\ldots,N\right\}  $, as well as $\lim_{t\rightarrow\infty}e^{-\left\vert
\xi\right\vert t}\psi\left(  t\right)  =0$ (by Lemma~\ref{le:pol-growth}), it
follows that $J_{\xi}=0$. \endproof\medskip

\noindent\textbf{Proof of Theorem \ref{thm:EU-natLspline}.} Note that any
function $s\in\mathcal{S}_{p,\xi}\left(  \tau\right)  $ is uniquely determined
as the extension on $\mathbb{R}$ of a function $s\in C^{2p-2}\left[
t_{0},t_{N}\right]  $ which is in $\mathrm{Ker}\mathcal{L}_{\xi}$ on any
subinterval $\left(  t_{j-1},t_{j}\right)  $, $j\in\left\{  1,\ldots
,N\right\}  $, and which satisfies the following endpoint conditions:%
\begin{equation}
\left\{
\begin{array}
[c]{c}%
\frac{d^{m}}{dt^{m}}\left(  \frac{d}{dt}-\left\vert \xi\right\vert \right)
^{p}s\left(  t\right)  |_{t=t_{0}}=0,\\
\frac{d^{m}}{dt^{m}}\left(  \frac{d}{dt}+\left\vert \xi\right\vert \right)
^{p}s\left(  t\right)  |_{t=t_{N}}=0,
\end{array}
\right.  \quad m\in\left\{  0,1,\ldots,p-2\right\}  .\label{eq:endpt}%
\end{equation}
The restriction of such a function $s$ to each subinterval $\left(
t_{j-1},t_{j}\right)  $, $j\in\left\{  1,\ldots,N\right\}  $, is determined by
$2p$ coefficients. Therefore imposing on $s$ the continuity conditions of
class $C^{2p-2}$ at each interior knot $t_{1},\ldots,t_{N-1}$, as well as the
endpoint conditions (\ref{eq:endpt}) and the interpolation conditions
(\ref{eq:1D-intconds}), we obtain a system of $\left(  2p-1\right)  \left(
N-1\right)  +2\left(  p-1\right)  +N+1=2pN$ linear equations for as many coefficients.

To establish the nature of this system, we assume zero interpolation data:
$y_{j}=0$, $j\in\left\{  0,1,\ldots,N\right\}  $, in which case the system
becomes homogeneous. Letting $s$ be determined by an arbitrary solution of
this homogeneous system and $\sigma=\psi:=s$ in (\ref{eq:FI-natLspline}), we
obtain $s\left(  t\right)  \in\mathrm{Ker}\left(  \frac{d}{dt}-\left\vert
\xi\right\vert \right)  ^{p}$ for $t\in\mathbb{R}$. Since $s\left(
t_{j}\right)  =0$, $j\in\left\{  0,1,\ldots,N\right\}  $, and $N+1\geq p$, we
deduce $s\equiv0$. It follows that the above homogeneous system has only the
trivial solution, which proves the conclusion of the theorem. \endproof
\medskip

\noindent\textbf{Proof of Theorem \ref{thm:var-p1D}.} Letting $\sigma:=s$ and
$\psi:=f-s$ in Theorem~\ref{thm:FI-natLspline},%
\[
\int_{-\infty}^{\infty}\left(  \frac{d}{dt}-\left\vert \xi\right\vert \right)
^{p}s\left(  t\right)  \left(  \frac{d}{dt}-\left\vert \xi\right\vert \right)
^{p}\overline{\left(  f\left(  t\right)  -s\left(  t\right)  \right)  }dt=0,
\]
hence%
\begin{align*}
0  & \leq\int_{-\infty}^{\infty}\left\vert \left(  \frac{d}{dt}-\left\vert
\xi\right\vert \right)  ^{p}\left(  f\left(  t\right)  -s\left(  t\right)
\right)  \right\vert ^{2}dt\\
& =\int_{-\infty}^{\infty}\left\vert \left(  \frac{d}{dt}-\left\vert
\xi\right\vert \right)  ^{p}f\left(  t\right)  \right\vert ^{2}dt-\int
_{-\infty}^{\infty}\left\vert \left(  \frac{d}{dt}-\left\vert \xi\right\vert
\right)  ^{p}s\left(  t\right)  \right\vert ^{2}dt.
\end{align*}
The inequality can become equality only if $f\left(  t\right)  -s\left(
t\right)  \in\mathrm{Ker}\left(  \frac{d}{dt}-\left\vert \xi\right\vert
\right)  ^{p}$ for $t\in\mathbb{R}$. Since $f-s$ vanishes at the knots
$t_{0},\ldots,t_{N}$, this implies $f\equiv s$. \endproof

\subsection{Proof of Lemma \ref{le:identity}}

Let $J_{\xi}$ denote the integral on the left-hand side of (\ref{eq:identity}%
). Expanding the two brackets inside this integral and collecting all the
terms that have the same power $l$ of $\left\vert \xi\right\vert $, we obtain%
\begin{equation}
J_{\xi}=\sum_{l=0}^{2p}\left(  -1\right)  ^{l}\left\vert \xi\right\vert
^{l}J_{\xi,l},\label{eq:I-xi-sum}%
\end{equation}
where%
\[
J_{\xi,l}:=\sum_{r=0}^{l}\left(
\begin{array}
[c]{c}%
p\\
l-r
\end{array}
\right)  \left(
\begin{array}
[c]{c}%
p\\
r
\end{array}
\right)  \int_{-\infty}^{\infty}\sigma^{\left(  p-r\right)  }\left(  t\right)
\,\overline{\psi^{\left(  p-l+r\right)  }\left(  t\right)  }\,dt.
\]
By convention, the two binomial coefficients in the last formula are zero if
$l-r>p$ or $r>p$, respectively. The lemma will follow by showing that in
(\ref{eq:I-xi-sum}) the odd power terms vanish and each of the even power
terms equals its correspondent term on the right-hand side of
(\ref{eq:identity}).

First, if $l=2k+1$, where $k\in\left\{  0,1,\ldots,p-1\right\}  $, then the
term of index $r$ of the sum $J_{\xi,l}$ contains the same product of binomial
coefficients as the term of index $l-r$. Further, for $r\leq k$, integration
by parts gives the relations%
\begin{align*}
& \int_{-\infty}^{\infty}\sigma^{\left(  p-r\right)  }\left(  t\right)
\,\overline{\psi^{\left(  p-l+r\right)  }\left(  t\right)  }\,dt\\
& =-\int_{-\infty}^{\infty}\sigma^{\left(  p-r-1\right)  }\left(  t\right)
\,\overline{\psi^{\left(  p-l+r+1\right)  }\left(  t\right)  }\,dt\\
& =\ldots=\left(  -1\right)  ^{l-2r}\int_{-\infty}^{\infty}\sigma^{\left(
p-l+r\right)  }\left(  t\right)  \,\overline{\psi^{\left(  p-r\right)
}\left(  t\right)  }\,dt.
\end{align*}
Indeed, the boundary terms in each integration by parts are zero due to the
exponential decay at $\pm\infty$ of $\sigma\left(  t\right)  $ and its
derivatives, as well as to the algebraic growth at $\pm\infty$ of $\psi\left(
t\right)  $ and its derivatives of order at most $p-1$. The latter growth
property is a consequence of Lemma~\ref{le:pol-growth} and the hypotheses on
$\psi$. Since $l-2r$ is odd, we deduce that the terms of indices $r$ and $l-r$
of the sum $J_{\xi,l}$ cancel, hence $J_{\xi,2k+1}=0$.

Second, if $l=2k$, where $k\in\left\{  0,1,\ldots,p\right\}  $, then, for
$0\leq r<k$, integrating by parts $k-r$ times yields%
\begin{align*}
& \int_{-\infty}^{\infty}\sigma^{\left(  p-r\right)  }\left(  t\right)
\,\overline{\psi^{\left(  p-2k+r\right)  }\left(  t\right)  }\,dt\\
& =-\int_{-\infty}^{\infty}\sigma^{\left(  p-r-1\right)  }\left(  t\right)
\,\overline{\psi^{\left(  p-2k+r+1\right)  }\left(  t\right)  }\,dt\\
& =\ldots=\left(  -1\right)  ^{k-r}\int_{-\infty}^{\infty}\sigma^{\left(
p-k\right)  }\left(  t\right)  \,\overline{\psi^{\left(  p-k\right)  }\left(
t\right)  }\,dt,
\end{align*}
where the boundary terms in each integration by parts are zero by the same
reasons as in the previous paragraph. Combining this with a similar argument
for $k<r\leq2k$, it follows that $J_{\xi,2k}$ is a multiple of the integral%
\[
\int_{-\infty}^{\infty}\sigma^{\left(  p-k\right)  }\left(  t\right)
\,\overline{\psi^{\left(  p-k\right)  }\left(  t\right)  }\,dt,
\]
the actual factor that multiplies the integral being the value%
\[
\left(  -1\right)  ^{k}\sum_{r=0}^{2k}\left(  -1\right)  ^{r}\left(
\begin{array}
[c]{c}%
p\\
2k-r
\end{array}
\right)  \left(
\begin{array}
[c]{c}%
p\\
r
\end{array}
\right)  =\left(
\begin{array}
[c]{c}%
p\\
k
\end{array}
\right)  .
\]
The last equality is a consequence of binomial expansions in the formal
identity $\left(  1-X\right)  ^{p}\left(  1+X\right)  ^{p}=\left(
1-X^{2}\right)  ^{p}$. Therefore (\ref{eq:identity}) is true. \endproof

\subsection{Proof of Theorem~\ref{thm:stability}}

For each $\xi\in\mathbb{R}^{n}\backslash\left\{  0\right\}  $, we first
introduce the integrable fundamental solution $\varphi_{\xi}$ of the operator
(\ref{eq:Lxi}), i.e. $\mathcal{L}_{\xi}\varphi_{\xi}=\delta_{0}$ (the Dirac
mass at the origin).

\begin{lemma}
The unique integrable fundamental solution of the operator $\mathcal{L}_{\xi}$
is%
\[
\varphi_{\xi}\left(  t\right)  =\frac{\left(  -1\right)  ^{p}}{\gamma
_{p}\left\vert \xi\right\vert ^{2p-1}}e^{-\left\vert \xi\right\vert \left\vert
t\right\vert }\sum_{l=0}^{p-1}c_{l}\left\vert \xi\right\vert ^{l}\left\vert
t\right\vert ^{l},\quad t\in\mathbb{R},
\]
where $c_{l}=\frac{\left(  2p-2-l\right)  !2^{l}}{l!\left(  p-1-l\right)  !}$,
$\forall l\in\left\{  0,\ldots,p-1\right\}  $, and $\gamma_{p}=\left(
p-1\right)  !2^{2p-1}$.
\end{lemma}

\noindent\proof Note that an integrable function $\varphi_{\xi}$ satisfies the
distributional equation $\mathcal{L}_{\xi}\varphi_{\xi}=\delta_{0}$ if and
only if its Fourier transform $\widehat{\varphi}_{\xi}\left(  u\right)
=\int_{\mathbb{R}}e^{-itu}\varphi_{\xi}\left(  t\right)  dt$ is%
\[
\widehat{\varphi}_{\xi}\left(  u\right)  =\frac{\left(  -1\right)  ^{p}%
}{\left(  u^{2}+\left\vert \xi\right\vert ^{2}\right)  ^{p}},\quad
u\in\mathbb{R}.
\]
By Fourier inversion,%
\begin{equation}
\varphi_{\xi}\left(  t\right)  =\frac{1}{2\pi}\int_{\mathbb{R}}e^{itu}%
\widehat{\varphi}_{\xi}\left(  u\right)  du=\frac{\left(  -1\right)  ^{p}}%
{\pi}\int_{0}^{\infty}\frac{\cos tu}{\left(  u^{2}+\left\vert \xi\right\vert
^{2}\right)  ^{p}}du.\label{eq:IFT}%
\end{equation}
The last integral is evaluated in Watson \cite[\S 6.16(1)]{Watson} as%
\[
\int_{0}^{\infty}\frac{\cos tu}{\left(  u^{2}+\left\vert \xi\right\vert
^{2}\right)  ^{p}}du=\frac{\pi^{1/2}}{\left(  p-1\right)  !}\left(
\frac{\left\vert t\right\vert }{2\left\vert \xi\right\vert }\right)
^{p-\frac{1}{2}}K_{p-\frac{1}{2}}\left(  \left\vert \xi\right\vert \left\vert
t\right\vert \right)  ,\quad t\not =0,
\]
where $K_{p-\frac{1}{2}}$ is the modified Bessel function expressible in
finite terms \cite[\S 3.71(12)]{Watson} by%
\[
K_{p-\frac{1}{2}}\left(  z\right)  =\left(  \frac{\pi}{2z}\right)
^{1/2}e^{-z}\sum_{l=0}^{p-1}\frac{\left(  p-1+l\right)  !}{l!\left(
p-1-l\right)  !\left(  2z\right)  ^{l}}.
\]
Therefore%
\[
\int_{0}^{\infty}\frac{\cos tu}{\left(  u^{2}+\left\vert \xi\right\vert
^{2}\right)  ^{p}}du=\frac{\pi e^{-\left\vert \xi\right\vert \left\vert
t\right\vert }}{\left(  p-1\right)  !\left(  2\left\vert \xi\right\vert
\right)  ^{2p-1}}\sum_{l=0}^{p-1}c_{l}\left\vert \xi\right\vert ^{l}\left\vert
t\right\vert ^{l},
\]
with the coefficients $c_{l}$ from the statement. The expression of
$\varphi_{\xi}$ follows.\endproof \medskip

For convenience, we will work with the following `normalization' of
$\varphi_{\xi}$:%
\begin{equation}
\widetilde{\varphi}_{\xi}\left(  t\right)  :=\left(  -1\right)  ^{p}\gamma
_{p}\left\vert \xi\right\vert ^{2p-1}\varphi_{\xi}\left(  t\right)
=e^{-\left\vert \xi\right\vert \left\vert t\right\vert }\sum_{l=0}^{p-1}%
c_{l}\left\vert \xi\right\vert ^{l}\left\vert t\right\vert ^{l},\quad
t\in\mathbb{R}.\label{eq:pd-RBF}%
\end{equation}
For each $\xi\in\mathbb{R}^{n}\backslash\left\{  0\right\}  $ and
$j\in\left\{  0,1,\ldots,N\right\}  $, define the function $L_{\xi,j}$
(independently of Theorem~\ref{thm:EU-natLspline}) by the `radial basis
function' (RBF) representation%
\begin{equation}
L_{\xi,j}\left(  t\right)  =\sum_{k=0}^{N}a_{jk}\widetilde{\varphi}_{\xi
}\left(  t-t_{k}\right)  ,\quad t\in\mathbb{R},\label{eq:Lagrange-RBF}%
\end{equation}
where the set of coefficients $\left\{  a_{jk}:k=0,\ldots,N\right\}  $ is
determined uniquely in the next lemma from the system of interpolation
conditions (\ref{eq:Lagrange-int}). Let%
\[
M_{\xi}:=\left(  \widetilde{\varphi}_{\xi}\left(  t_{j}-t_{k}\right)  \right)
_{j,k=0}^{N}%
\]
be the matrix of this interpolation system under representation
(\ref{eq:Lagrange-RBF}).

\begin{lemma}
For every $\xi\in\mathbb{R}^{n}\backslash\left\{  0\right\}  $, the matrix
$M_{\xi}$ is positive definite, hence nonsingular. Accordingly, for each
$j\in\left\{  0,\ldots,N\right\}  $, the function $L_{\xi,j}$ defined by
(\ref{eq:Lagrange-RBF}) and (\ref{eq:Lagrange-int}) is the Lagrange function
in $\mathcal{S}_{p,\xi}\left(  \tau\right)  $ arising from
Theorem~\ref{thm:EU-natLspline}.
\end{lemma}

\noindent\proof Writing the Fourier inversion formula (\ref{eq:IFT})\ in the
form%
\[
\widetilde{\varphi}_{\xi}\left(  t\right)  =\frac{\gamma_{p}\left\vert
\xi\right\vert ^{2p-1}}{2\pi}\int_{\mathbb{R}}\frac{e^{itu}}{\left(
u^{2}+\left\vert \xi\right\vert ^{2}\right)  ^{p}}du,\quad t\in\mathbb{R},
\]
we obtain, for any column vector $\mathbf{v}=\left(  v_{0},\ldots
,v_{N}\right)  ^{T}\in\mathbb{R}^{N+1}$,%
\begin{align*}
\mathbf{v}^{T}M_{\xi}\mathbf{v}  & =\sum_{j=0}^{N}\sum_{k=0}^{N}%
\widetilde{\varphi}_{\xi}\left(  t_{j}-t_{k}\right)  v_{j}v_{k}\\
& =\frac{\gamma_{p}\left\vert \xi\right\vert ^{2p-1}}{2\pi}\int_{\mathbb{R}%
}\frac{\left\vert \sum_{k=0}^{N}v_{k}e^{it_{k}u}\right\vert ^{2}}{\left(
u^{2}+\left\vert \xi\right\vert ^{2}\right)  ^{p}}du\geq0.
\end{align*}
Further, since the functions $e^{it_{k}u}$, for $k\in\left\{  0,\ldots
,N\right\}  $, are linearly independent (their Wronskian being the multiple of
a Vandermonde determinant), it follows that the above inequality is strict if
$\mathbf{v}\not =0$. Therefore the matrix $M_{\xi}$ is positive definite and nonsingular.

For the last part of the lemma, it is sufficient to prove that $L_{\xi,j}$
defined by (\ref{eq:Lagrange-RBF}) and (\ref{eq:Lagrange-int}) is in
$\mathcal{S}_{p,\xi}\left(  \tau\right)  $. Differentiating the above Fourier
inversion formula shows that $\widetilde{\varphi}_{\xi}\in C^{2p-2}\left(
\mathbb{R}\right)  $, so the linear combination (\ref{eq:Lagrange-RBF}) is in
$C^{2p-2}\left(  \mathbb{R}\right)  $. On the other hand, for any non-negative
integer $l$, any knot $t_{k}\in\tau$, and any parameter $\xi$, the following
formula shows that the translate of an exponential polynomial is spanned by
exponential polynomials of the same type, namely%
\[
\left(  t-t_{k}\right)  ^{l}e^{\pm\left\vert \xi\right\vert \left(
t-t_{k}\right)  }=\sum_{r=0}^{l}a_{r}\left(  l,t_{k},\left\vert \xi\right\vert
\right)  t^{r}e^{\pm\left\vert \xi\right\vert t},
\]
for some coefficients $a_{r}=a_{r}\left(  l,t_{k},\left\vert \xi\right\vert
\right)  $ (cf.\ \cite{UB03}). We deduce that the form (\ref{eq:Lagrange-RBF})
of $L_{\xi,j}$ also satisfies conditions \textit{(ii)} and \textit{(iii)} of
Definition~\ref{def:natLspline}, as required. The lemma is proved. \endproof

\begin{lemma}
\label{le:diag-dom}There exists a constant $\mu=\mu\left(  p,\tau\right)  >0 $
such that, for any $\xi\in\mathbb{R}^{n}$ with $\left\vert \xi\right\vert
\geq\mu$ and any eigenvalue $\lambda$ of $M_{\xi}$, we have $\lambda\geq
\frac{1}{2}c_{0}$, where $c_{0}=\frac{\left(  2p-2\right)  !}{\left(
p-1\right)  !}$.
\end{lemma}

\noindent\proof We follow an idea from Narcowich-Ward \cite{NW91} and Schaback
\cite{RS95}. Note that the diagonal entries of $M_{\xi}$ are all equal to
$\widetilde{\varphi}_{\xi}\left(  0\right)  =c_{0}$. We show that $M_{\xi}$ is
diagonally dominant for sufficiently large $\left\vert \xi\right\vert $. Let%
\[
\psi\left(  u\right)  :=e^{-u}\sum_{l=0}^{p-1}c_{l}u^{l},\quad u\geq0,
\]
so that $\widetilde{\varphi}_{\xi}\left(  t\right)  =\psi\left(  \left\vert
\xi\right\vert \left\vert t\right\vert \right)  $, $\forall t\in\mathbb{R}$.
Since $\lim_{u\rightarrow\infty}\psi\left(  u\right)  =0$, there exists
$\delta>0$ depending only on $p$ and $N$ such that%
\[
\psi\left(  u\right)  \leq\frac{c_{0}}{2N},\quad\forall u\geq\delta.
\]
Let $\mu:=\delta/ \min\limits_{j\not =k}\left\vert t_{j}-t_{k}\right\vert $.
Then for all $\xi$ with $\left\vert \xi\right\vert \geq\mu$ and all $j\not =%
k$, we have $\left\vert \xi\right\vert \left\vert t_{j}-t_{k}\right\vert
\geq\delta$, from which%
\[
\rho:=\max_{0\leq j\leq N}\sum_{\substack{k=0 \\k\not =j}}^{N}\widetilde
{\varphi}_{\xi}\left(  t_{j}-t_{k}\right)  \leq\frac{c_{0}}{2}.
\]
It follows that the matrix $M_{\xi}$ is diagonally dominant and Gershgorin's
circle theorem \cite[p.~395]{GVL96} implies $\lambda\geq\widetilde{\varphi
}_{\xi}\left(  0\right)  -\rho\geq\frac{1}{2}c_{0}$, for any eigenvalue
$\lambda$ of $M_{\xi}$. \endproof \medskip

\noindent\textbf{Proof of Theorem~\ref{thm:stability}.} First we show that
there exists a constant $A=A\left(  p,\tau\right)  $ such that the
coefficients $a_{jk}$ of the RBF representation (\ref{eq:Lagrange-RBF})
satisfy%
\begin{equation}
\left\vert a_{jk}\right\vert \leq A, \quad\forall j,k\in\left\{
0,\ldots,N\right\}  ,\ \forall\left\vert \xi\right\vert \geq
1/2.\label{eq:unif-bd}%
\end{equation}
Note that, if $\mathbf{a}_{j}:=\left(  a_{j0},\ldots,a_{jN}\right)  ^{T}%
\in\mathbb{R}^{N+1}$ and $\mathbf{e}_{j}$ is the $j$th column of the identity
matrix of order $N+1$, then $\mathbf{a}_{j}=M_{\xi}^{-1}\mathbf{e}_{j}$, i.e.
$\mathbf{a}_{j}$ is the $j$th column of the inverse matrix $M_{\xi}^{-1}$.
Hence, for $\left\vert \xi\right\vert \geq\mu$ and any $j$, $k$,
Lemma~\ref{le:diag-dom} implies%
\[
\left\vert a_{jk}\right\vert \leq\left\Vert \mathbf{a}_{j}\right\Vert
_{2}=\left\Vert M_{\xi}^{-1}\mathbf{e}_{j}\right\Vert _{2}\leq\left\Vert
M_{\xi}^{-1}\right\Vert _{2}=\frac{1}{\min\left\vert \lambda\right\vert }%
\leq\frac{2}{c_{0}},
\]
where $\left\Vert \mathbf{a}_{j}\right\Vert _{2}$ is the Euclidean norm in
$\mathbb{R}^{N+1}$, $\left\Vert M_{\xi}^{-1}\right\Vert _{2}$ is the induced
matrix norm, and the minimum is taken over all eigenvalues $\lambda$ of
$M_{\xi}$. If $\frac{1}{2}\geq\mu$, then (\ref{eq:unif-bd}) holds with
$A:=\frac{2}{c_{0}}$. If $\frac{1}{2}<\mu$, then for each $\xi\in
\mathbb{R}^{n}$ with $\frac{1}{2}\leq\left\vert \xi\right\vert \leq\mu$, we
use Cramer's rule to express each coefficient $a_{jk}$ as the ratio of two
determinants, the denominator being $\det M_{\xi}$. The numerator is the
determinant of the matrix obtained from $M_{\xi}$ by replacing its $j$th
column with $\mathbf{e}_{k}$. Since the entries $\widetilde{\varphi}_{\xi
}\left(  t_{j}-t_{k}\right)  $ of $M_{\xi}$ depend continuously on $\xi$, it
follows that there exists a constant $b_{0}=b_{0}\left(  p,\tau\right)  >0$,
such that $\left\vert a_{jk}\right\vert \leq b_{0}$ for all $j$, $k$, and all
$\xi$ with $\frac{1}{2}\leq\left\vert \xi\right\vert \leq\mu$. Therefore
(\ref{eq:unif-bd}) is obtained by letting $A:=\max\left\{  b_{0},\frac
{2}{c_{0}}\right\}  $.

Since the values $\widetilde{\varphi}_{\xi}\left(  t\right)  $ of
(\ref{eq:pd-RBF}) are positive and bounded above by a constant for all
$t\in\mathbb{R}$ and $\xi\in\mathbb{R}^{n}\backslash\left\{  0\right\}  $,
from (\ref{eq:Lagrange-RBF}) and (\ref{eq:unif-bd}) we deduce that
(\ref{eq:stability}) holds for $m=0$. The corresponding estimates for the
derivatives of $L_{\xi,j}$ follow in the same way, since (\ref{eq:pd-RBF})
implies%
\begin{equation}
\widetilde{\varphi}_{\xi}^{\left(  m\right)  }\left(  t\right)  =\left\vert
\xi\right\vert ^{m}e^{-\left\vert \xi\right\vert \left\vert t\right\vert }%
\sum_{l=0}^{p-1}c_{m,l}\left\vert \xi\right\vert ^{l}\left\vert t\right\vert
^{l},\quad t\in\mathbb{R},\ m\in\left\{  1,\ldots,2p-2\right\}
,\label{eq:RBF-der}%
\end{equation}
for certain coefficients $c_{m,l}$. This completes the proof. \endproof

\section{Transfinite interpolation with Beppo Levi poly\-splines}

\label{polysplines}

For any non-negative integer $r$, let $W^{r}\left(  \mathbb{T}^{n}\right)  $
be the space of all continuous functions $f:\mathbb{R}^{n}\rightarrow
\mathbb{C}$ which are $2\pi$-periodic in each variable and satisfy the
condition%
\[
\left\Vert f\right\Vert _{r}:=\sum_{\xi\in\mathbb{Z}^{n}}\left\vert
\widehat{f}_{\xi}\right\vert \left(  1+\left\vert \xi\right\vert \right)
^{r}<\infty,
\]
where $\widehat{f}_{\xi}:=\int_{\mathbb{T}^{n}}e^{-i\left\langle
\xi,y\right\rangle }f\left(  y\right)  dy$, for $\xi\in\mathbb{Z}^{n}$, are
the Fourier coefficients of $f$. In particular, $W^{0}\left(  \mathbb{T}%
^{n}\right)  $ is the Wiener algebra of functions with absolutely convergent
Fourier series. Note that $0\leq r_{1}\leq r_{2}$ implies $\left\Vert
f\right\Vert _{r_{1}}\leq\left\Vert f\right\Vert _{r_{2}}$, hence $W^{r_{2}%
}\left(  \mathbb{T}^{n}\right)  \subset W^{r_{1}}\left(  \mathbb{T}%
^{n}\right)  $. It is also straightforward that $C^{k}\left(  \mathbb{T}%
^{n}\right)  \subset W^{r}\left(  \mathbb{T}^{n}\right)  \subset C^{r}\left(
\mathbb{T}^{n}\right)  $, where $k$ is the least integer greater than or equal
to $r+\frac{n+1}{2}$.

Now recall the setting of the Introduction, in which $p\geq2$. First, we state
the existence and uniqueness of transfinite interpolation with Beppo Levi
polysplines of order $p$ on strips.

\begin{theorem}
\label{thm:EU-poly-nat}For any set of data functions $f_{j}\in W^{2p-2}\left(
\mathbb{T}^{n}\right)  $, $j\in\left\{  0,\ldots,N\right\}  $, there exists a
unique Beppo Levi polyspline\ $S\in\mathcal{S}_{p}\left(  \tau,n\right)  $
such that%
\begin{equation}
S\left(  t_{j},y\right)  =f_{j}\left(  y\right)  ,\quad j\in\left\{
0,1,\ldots,N\right\}  ,\ y\in\mathbb{T}^{n}.\label{eq:transfinite-int}%
\end{equation}

\end{theorem}

Let $\mathcal{B}_{p}:=\mathcal{B}_{p}\left(  \mathbb{R}\times\mathbb{T}%
^{n}\right)  $ be the space of all functions $F\in C^{p}\left(  \mathbb{R}%
^{n+1}\right)  $ that are $2\pi$-periodic in each of the last $n$ variables
and satisfy, in usual multi-index notation,%
\[
\partial^{\alpha}F\in L^{2}\left(  \mathbb{R}\times\mathbb{T}^{n}\right)
,\quad\forall\left\vert \alpha\right\vert =p.
\]
On $\mathcal{B}_{p}$ we define the following semi-inner product and induced
seminorm:%
\begin{align}
\left\langle F,G\right\rangle _{\mathcal{B}_{p}} &  :=\int\limits_{\mathbb{R}%
\times\mathbb{T}^{n}}\sum_{\left\vert \alpha\right\vert =p}\frac{p!}{\alpha
!}\partial^{\alpha}F\,\overline{\partial^{\alpha}G}\,dtdy,\quad\forall
F,G\in\mathcal{B}_{p},\label{eq:Duchon}\\
\left\Vert F\right\Vert _{\mathcal{B}_{p}} &  :=\left\langle F,F\right\rangle
_{\mathcal{B}_{p}}^{1/2},\quad\forall F\in\mathcal{B}_{p}.\nonumber
\end{align}
Note that $\mathcal{S}_{p}\left(  \tau,n\right)  \subset\mathcal{B}_{p}$ for
$p\geq2$, and that $\left\Vert F\right\Vert _{\mathcal{B}_{p}}^{2}$ coincides
with (\ref{eq:Du-per}).

The following orthogonality result is the polyspline analog of
Theorem~\ref{thm:FI-natLspline}.

\begin{theorem}
\label{thm:poly-FI}Let $S\in\mathcal{S}_{p}\left(  \tau,n\right)  $ be any
Beppo Levi polyspline as in Definition~\ref{def:poly-nat}. If $G\in
\mathcal{B}_{p}$ satisfies $G\left(  t_{j},y\right)  =0$, $\forall
j\in\left\{  0,1,\ldots,N\right\}  $, $\forall y\in\mathbb{T}^{n}$, then%
\[
\left\langle S,G\right\rangle _{\mathcal{B}_{p}}=0.
\]

\end{theorem}

\noindent This `fundamental identity' is eventually used to establish the
characterization of the Beppo Levi polyspline $S$ of
Theorem~\ref{thm:EU-poly-nat} as minimizer of the Duchon type seminorm
$\left\Vert \cdot\right\Vert _{\mathcal{B}_{p}}$ subject to the transfinite
interpolation conditions (\ref{eq:transfinite-int}).

\begin{theorem}
\label{thm:var-p}Given the set of data functions $\left\{  f_{0},f_{1}%
,\ldots,f_{N}\right\}  \subset W^{2p-2}\left(  \mathbb{T}^{n}\right)  $, let
$S\in\mathcal{S}_{p}\left(  \tau\times\mathbb{T}^{n}\right)  $ be the Beppo
Levi polyspline that satisfies conditions (\ref{eq:transfinite-int}). If
$F\in\mathcal{B}_{p}$ is any other function satisfying $F\left(
t_{j},y\right)  =f_{j}\left(  y\right)  $, $\forall j\in\left\{
0,\ldots,N\right\}  $,$\ \forall y\in\mathbb{T}^{n}$, then%
\[
\left\Vert S\right\Vert _{\mathcal{B}_{p}} <\left\Vert F\right\Vert
_{\mathcal{B}_{p}}.
\]

\end{theorem}

The proofs below will employ standard notation. In particular, $\mathbb{Z}%
_{+}:=\left\{  0,1,\ldots\right\}  $ and, for a multi-index $\beta=\left(
\beta_{1},\ldots,\beta_{n}\right)  \in\mathbb{Z}_{+}^{n}$, we let $\left\vert
\beta\right\vert :=\beta_{1}+\ldots+\beta_{n}$, $\beta! := \beta_{1} !
\ldots\beta_{n} ! $, and $\partial_{y}^{\beta}:=\frac{\partial^{\beta_{1}}%
}{\partial y_{1}^{\beta_{1}}}\ldots\frac{\partial^{\beta_{n}}}{\partial
y_{n}^{\beta_{n}}}$.

\subsection{Proof of Theorem~\ref{thm:EU-poly-nat}}

To construct the natural polyspline $S$ with the properties stated in
Theorem~\ref{thm:EU-poly-nat}, we start from the absolutely convergent Fourier
series representations%
\[
f_{j}\left(  y\right)  =\sum_{\xi\in\mathbb{Z}^{n}}\widehat{f}_{j,\xi
}e^{i\left\langle \xi,y\right\rangle },\quad j\in\left\{  0,1,\ldots
,N\right\}  ,\ y\in\mathbb{T}^{n},
\]
where $\widehat{f}_{j,\xi}$, for $\xi\in\mathbb{Z}^{n}$, are the Fourier
coefficients of $f_{j}\in W^{2p-2}\left(  \mathbb{T}^{n}\right)  $. For each
$\xi\in\mathbb{Z}^{n}$, by Theorem~\ref{thm:EU-natLspline} and its classical
version for $\xi=0$, let $\widehat{S}_{\xi}\in\mathcal{S}_{p,\xi}\left(
\tau\right)  $ be the unique natural $\mathcal{L}_{\xi}$-spline on $\tau$,
such that%
\[
\widehat{S}_{\xi}\left(  t_{j}\right)  =\widehat{f}_{j,\xi},\quad j\in\left\{
0,1,\ldots,N\right\}  .
\]
By the Lagrange formula (\ref{eq:Lagrange}), we have%
\begin{equation}
\widehat{S}_{\xi}\left(  t\right)  =\sum_{j=0}^{N}\widehat{f}_{j,\xi}L_{\xi
,j}\left(  t\right)  ,\quad t\in\mathbb{R}.\label{eq:Lag-rep}%
\end{equation}
We will prove that the function $S$ defined by%
\begin{equation}
S\left(  t,y\right)  :=\sum_{\xi\in\mathbb{Z}^{n}}\widehat{S}_{\xi}\left(
t\right)  e^{i\left\langle \xi,y\right\rangle },\quad y\in\mathbb{T}%
^{n},\ t\in\mathbb{R},\label{eq:F-series}%
\end{equation}
possesses all the properties required by Theorem~\ref{thm:EU-poly-nat}.

First, we show that $S\in C^{2p-2}\left(  \mathbb{R}\times\mathbb{T}%
^{n}\right)  $ and $S$ satisfies (\ref{eq:transfinite-int}). Let
$m\in\mathbb{Z}_{+}$, $\beta\in\mathbb{Z}_{+}^{n}$, and $\alpha=\left(
m,\beta\right)  $, such that $\left\vert \alpha\right\vert =m+\left\vert
\beta\right\vert \leq2p-2$. For $t\in\left[  t_{0},t_{N}\right]  $, we
substitute (\ref{eq:Lag-rep}) into (\ref{eq:F-series})\ and use the estimate
(\ref{eq:stability}) to obtain formally%
\begin{align*}
\left\vert \partial^{\alpha}S\left(  t,y\right)  \right\vert  &  \leq\sum
_{\xi\in\mathbb{Z}^{n}}\sum_{j=0}^{N}\left\vert \widehat{f}_{j,\xi}\right\vert
\left\vert \partial_{t}^{m}\partial_{y}^{\beta}L_{\xi,j}\left(  t\right)
e^{i\left\langle \xi,y\right\rangle }\right\vert \\
&  \leq C_{0}\sum_{\xi\in\mathbb{Z}^{n}}\sum_{j=0}^{N}\left\vert \widehat
{f}_{j,\xi}\right\vert \left(  1+\left\vert \xi\right\vert ^{m}\right)
\left\vert \xi\right\vert ^{\left\vert \beta\right\vert }\\
&  \leq2C_{0}\sum_{j=0}^{N}\left\Vert f_{j}\right\Vert _{\left\vert
\alpha\right\vert },\quad\forall y\in\mathbb{T}^{n}.
\end{align*}
For $t\not \in \left[  t_{0},t_{N}\right]  $, we use the same estimates for
all terms corresponding to $\xi\in\mathbb{Z}^{n}\backslash\left\{  0\right\}
$ in the first line of the above display, while the modulus of the term
corresponding to $\xi=0$ is bounded above by%
\[
\sum_{j=0}^{N}\left\vert \widehat{f}_{j,0}\right\vert \left\vert \partial
_{t}^{m}\partial_{y}^{\beta}L_{0,j}\left(  t\right)  \right\vert .
\]
This expression vanishes if $\left\vert \beta\right\vert \geq1$, or if
$\beta=0$ and $m\geq p$, since $L_{0,j}$ is a polynomial of degree at most
$p-1$ outside the interval $\left[  t_{0},t_{N}\right]  $. In the remaining
case $\beta=0$ and $0\leq m\leq p-1$, we use the polynomial growth estimate%
\[
\left\vert \frac{d^{m}}{dt^{m}}L_{0,j}\left(  t\right)  \right\vert \leq
C_{1}\left(  1+\left\vert t\right\vert ^{p-1-m}\right)  ,\quad\forall
t\not \in \left[  t_{0},t_{N}\right]  ,\,\forall j\in\left\{  0,\ldots
,N\right\}  ,
\]
with the constant $C_{1}$ depending only on $p$ and $\tau$. Altogether, these
estimates show that the series (\ref{eq:F-series}) is absolutely and uniformly
convergent on compact sets in $\mathbb{R}\times\mathbb{T}^{n}$ and can be
differentiated termwise up to the total order $2p-2$. We also obtain
(\ref{eq:transfinite-int}) from the Lagrange conditions (\ref{eq:Lagrange-int}%
) satisfied by $L_{\xi,j}$.

Next, to show that $S$ is polyharmonic of order $p$ on each of the strips
$\Omega_{j}$, $j\in\left\{  0,\ldots,N+1\right\}  $, we use the following result.

\begin{lemma}
\label{le:nicolescu}Suppose that a series of functions is absolutely and
uniformly convergent on compact sets in the open domain $\Omega$ and each term
of the series is polyharmonic of order $p$ on $\Omega$. Then the sum of the
series is also polyharmonic of order $p$ on $\Omega$.
\end{lemma}

\noindent This is well-known in the harmonic case $p=1$. For $p\geq2$, it was
proved by Nicolesco \cite[p.~23]{Ni36} based on integral mean representations
of polyharmonic functions. In our case, each term of (\ref{eq:F-series})
satisfies%
\begin{align*}
\Delta^{p}\left[  \widehat{S}_{\xi}\left(  t\right)  e^{i\left\langle
\xi,y\right\rangle }\right]   &  =\left(  \frac{\partial^{2}}{\partial t^{2}%
}+\sum_{\nu=1}^{d}\frac{\partial^{2}}{\partial y_{\nu}^{2}}\right)
^{p}\left[  \widehat{S}_{\xi}\left(  t\right)  e^{i\left\langle \xi
,y\right\rangle }\right] \\
&  =e^{i\left\langle \xi,y\right\rangle }\left(  \frac{d^{2}}{dt^{2}%
}-\left\vert \xi\right\vert ^{2}\right)  ^{p}\widehat{S}_{\xi}\left(
t\right)  =0,
\end{align*}
for $\left(  t,y\right)  \in\Omega_{j}$, $j\in\left\{  0,\ldots,N+1\right\}
$, since $\widehat{S}_{\xi}\in\mathcal{S}_{p,\xi}\left(  \tau\right)  $.
Therefore $S$ is piecewise polyharmonic as required.

It remains to prove that $S$ satisfies the Beppo Levi condition
(\ref{eq:beppolevi}). Let $m\in\mathbb{Z}_{+}$, $\beta\in\mathbb{Z}_{+}^{n}$,
and $\alpha=\left(  m,\beta\right)  $, such that $\left\vert \alpha\right\vert
=m+\left\vert \beta\right\vert =p$. It is sufficient to prove that the
termwise partial derivative $\partial^{\alpha}$ of series (\ref{eq:F-series})
is a Cauchy series in $L^{2}\left(  \left(  \mathbb{R}\backslash\left[
t_{0},t_{N}\right]  \right)  \times\mathbb{T}^{n}\right)  $, which implies
convergence to its sum $\partial^{\alpha}S\left(  t,y\right)  $ in
$L^{2}\left(  \left(  \mathbb{R}\backslash\left[  t_{0},t_{N}\right]  \right)
\times\mathbb{T}^{n}\right)  $. Note that the partial derivative
$\partial^{\alpha}$ of order $p$ of the term corresponding to $\xi=0$ in
(\ref{eq:F-series}) vanishes due to the natural conditions satisfied by
$\widehat{S}_{\xi}$ outside the interval $\left[  t_{0},t_{N}\right]  $. If
$\xi\in\mathbb{Z}^{n}\backslash\left\{  0\right\}  $ and $0\leq j\leq N$,
relations (\ref{eq:pd-RBF}) and (\ref{eq:Lagrange-RBF}) give the following
representation of $L_{\xi,j}$ for $t\leq t_{0}$:%
\[
L_{\xi,j}\left(  t\right)  =\sum_{k=0}^{N}a_{jk}e^{\left\vert \xi\right\vert
\left(  t-t_{k}\right)  }\sum_{l=0}^{p-1}c_{l}\left\vert \xi\right\vert
^{l}\left(  t_{k}-t\right)  ^{l}.
\]
Hence, for $m=0$ and $\left\vert \beta\right\vert =p$, Minkowski's inequality
and (\ref{eq:unif-bd}) imply%
\begin{align*}
& \left(  \int_{-\infty}^{t_{0}}\int_{\mathbb{T}^{n}}\left\vert \partial
_{y}^{\beta}L_{\xi,j}\left(  t\right)  e^{i\left\langle \xi,y\right\rangle
}\right\vert ^{2}dydt\right)  ^{1/2}\\
& \leq\left\vert \xi\right\vert ^{p}\sum_{k=0}^{N}\left\vert a_{jk}\right\vert
\sum_{l=0}^{p-1}c_{l}\left\vert \xi\right\vert ^{l}\left(  \int_{-\infty
}^{t_{0}}\left(  t_{k}-t\right)  ^{2l}e^{2\left\vert \xi\right\vert \left(
t-t_{k}\right)  }dt\right)  ^{1/2}\\
& \leq A\left\vert \xi\right\vert ^{p}\sum_{k=0}^{N}\sum_{l=0}^{p-1}%
c_{l}\left\vert \xi\right\vert ^{l}\left(  \int_{-\infty}^{t_{k}}\left(
t_{k}-t\right)  ^{2l}e^{2\left\vert \xi\right\vert \left(  t-t_{k}\right)
}dt\right)  ^{1/2}.
\end{align*}
Since, by changing the variable, the last integral becomes Euler's integral%
\[
\int_{0}^{\infty}u^{2l}e^{-2\left\vert \xi\right\vert u}du=\frac{\left(
2l\right)  !}{\left(  2\left\vert \xi\right\vert \right)  ^{2l+1}},
\]
we obtain%
\[
\left(  \int_{-\infty}^{t_{0}}\int_{\mathbb{T}^{n}}\left\vert \partial
_{y}^{\beta}L_{\xi,j}\left(  t\right)  e^{i\left\langle \xi,y\right\rangle
}\right\vert ^{2}dydt\right)  ^{1/2}\leq C_{2}\left\vert \xi\right\vert
^{p-\frac{1}{2}},
\]
for a constant $C_{2}$ depending only on $p$ and $\tau$. For $m\geq1$ and
$m+\left\vert \beta\right\vert =p$, a similar argument based on the RBF
derivative (\ref{eq:RBF-der}) gives%
\[
\left(  \int_{-\infty}^{t_{0}}\int_{\mathbb{T}^{n}}\left\vert \partial_{t}%
^{m}\partial_{y}^{\beta}L_{\xi,j}\left(  t\right)  e^{i\left\langle
\xi,y\right\rangle }\right\vert ^{2}dydt\right)  ^{1/2}\leq C_{2}\left\vert
\xi\right\vert ^{p-\frac{1}{2}},
\]
after increasing $C_{2}$ if necessary. Therefore, if $\alpha\in\mathbb{Z}%
_{+}^{n+1}$ is any multi-index with $\left\vert \alpha\right\vert =p$ and
$\mathcal{F}$ is any finite subset of $\mathbb{Z}^{n}\backslash\left\{
0\right\}  $, another application of Minkowski's inequality gives%
\begin{align*}
&  \left(  \int_{-\infty}^{t_{0}}\int_{\mathbb{T}^{n}}\left\vert
\sum\nolimits_{\xi\in\mathcal{F}}\sum_{j=0}^{n}\widehat{f}_{j,\xi}%
\partial^{\alpha}L_{\xi,j}\left(  t\right)  e^{i\left\langle \xi
,y\right\rangle }\right\vert ^{2}dydt\right)  ^{1/2}\\
&  \leq\sum\nolimits_{\xi\in\mathcal{F}}\sum_{j=0}^{n}\left\vert \widehat
{f}_{j,\xi}\right\vert \left(  \int_{-\infty}^{t_{0}}\int_{\mathbb{T}^{n}%
}\left\vert \partial^{\alpha}L_{\xi,j}\left(  t\right)  e^{i\left\langle
\xi,y\right\rangle }\right\vert ^{2}dydt\right)  ^{1/2}\\
&  \leq C_{2}\sum\nolimits_{\xi\in\mathcal{F}}\sum_{j=0}^{n}\left\vert
\widehat{f}_{j,\xi}\right\vert \left\vert \xi\right\vert ^{p-\frac{1}{2}}.
\end{align*}
Since $f_{j}\in W^{2p-2}\left(  \mathbb{T}^{n}\right)  $ and $p-\frac{1}%
{2}<2p-2$ for $2\leq p$, we deduce that any termwise partial derivative of
total order $p$ of the series (\ref{eq:F-series}) is a Cauchy series in
$L^{2}\left(  \left(  -\infty,t_{0}\right)  \times\mathbb{T}^{n}\right)  $. By
similar arguments, this conclusion also holds for the interval $\left(
t_{N},\infty\right)  $ in place of $\left(  -\infty,t_{0}\right)  $, which
completes the proof that $S$ satisfies the Beppo-Levi condition
(\ref{eq:beppolevi}).

The uniqueness of a natural polyspline $S$ satisfying the conditions of
Theorem~\ref{thm:EU-poly-nat} follows directly from the uniqueness of natural
$\mathcal{L}_{\xi}$-splines established in Theorem~\ref{thm:EU-natLspline}, as
in the proof of uniqueness of \cite[Theorem~15]{AB08}. \endproof

\subsection{Proofs of Theorems~\ref{thm:poly-FI} and~\ref{thm:var-p}}

The proof of Theorem~\ref{thm:poly-FI} will employ Parseval type
representations of the semi-inner product (\ref{eq:Duchon}) with respect to
the $y$ variable. For any $F\in\mathcal{B}_{p}$ and $t\in\mathbb{R}$, define
the Fourier coefficients of $F$ with respect to its last $n$ variables by%
\begin{equation}
\widehat{F}_{\xi}\left(  t\right)  :=\int_{\mathbb{T}^{n}}e^{-i\left\langle
\xi,y\right\rangle }F\left(  t,y\right)  \,dy,\quad\xi\in\mathbb{Z}%
^{n}.\label{eq:F-coeffs}%
\end{equation}

\begin{lemma}
\label{le:Parseval-Fubini}Let $F,G\in\mathcal{B}_{p}$. If $m\in\mathbb{Z}_{+}%
$, $\beta\in\mathbb{Z}_{+}^{n}$, and $\alpha=\left(  m,\beta\right)  $ satisfy
$\left\vert \alpha\right\vert =m+\left\vert \beta\right\vert =p$, then the
following identities hold with absolutely convergent integrals and series:%
\begin{align}
\int\limits_{\mathbb{T}^{n}}\left\vert \partial^{\alpha}F\left(  t,y\right)
\right\vert ^{2}dy  & =\sum_{\xi\in\mathbb{Z}^{n}}\xi^{2\beta}\,\left\vert
\frac{d^{m}}{dt^{m}}\widehat{F}_{\xi}\left(  t\right)  \right\vert
^{2},\label{eq:P1}\\
\int\limits_{\mathbb{T}^{n}}\partial^{\alpha}F\left(  t,y\right)
\,\overline{\partial^{\alpha}G\left(  t,y\right)  }\,dy  & =\sum_{\xi
\in\mathbb{Z}^{n}}\xi^{2\beta}\frac{d^{m}}{dt^{m}}\widehat{F}_{\xi}\left(
t\right)  \,\frac{d^{m}}{dt^{m}}\overline{\widehat{G}_{\xi}\left(  t\right)
},\label{eq:P2}\\
\int\limits_{\mathbb{R}}\int\limits_{\mathbb{T}^{n}}\left\vert \partial
^{\alpha}F\left(  t,y\right)  \right\vert ^{2}dydt  & =\sum_{\xi\in
\mathbb{Z}^{n}}\xi^{2\beta}\int\limits_{\mathbb{R}}\left\vert \frac{d^{m}%
}{dt^{m}}\widehat{F}_{\xi}\left(  t\right)  \right\vert ^{2}dt,\label{eq:P3}\\
\int\limits_{\mathbb{R}} \! \int\limits_{\mathbb{T}^{n}}\partial^{\alpha
}F\left(  t,y\right)  \overline{\partial^{\alpha}G\left(  t,y\right)  }dydt  &
=\sum_{\xi\in\mathbb{Z}^{n}}\xi^{2\beta} \! \int\limits_{\mathbb{R}}%
\frac{d^{m}}{dt^{m}}\widehat{F}_{\xi}\left(  t\right)  \frac{d^{m}}{dt^{m}%
}\overline{\widehat{G}_{\xi}\left(  t\right)  } dt.\label{eq:P4}%
\end{align}

\end{lemma}

\noindent\proof Differentiating with respect to $t$ in (\ref{eq:F-coeffs}),
then integrating by parts with respect to each component of $y$ and using the
periodicity of $F$ in $y$ and the hypothesis $F\in C^{p}\left(  \mathbb{R}%
^{n+1}\right)  $, we obtain%
\[
\left(  i\xi\right)  ^{\beta}\frac{d^{m}}{dt^{m}}\widehat{F}_{\xi}\left(
t\right)  =\int_{\mathbb{T}^{n}}e^{-i\left\langle \xi,y\right\rangle }%
\partial_{y}^{\beta}\partial_{t}^{m}F\left(  t,y\right)  dy,\quad\xi
\in\mathbb{Z}^{n},
\]
where $\xi\not =0$ if $\beta\not =0$. Since $\partial_{y}^{\beta}\partial
_{t}^{m}F\left(  t,\cdot\right)  \in C\left(  \mathbb{T}^{n}\right)  \subset
L^{2}\left(  \mathbb{T}^{n}\right)  $, the sequence of Fourier coefficients of
$\partial_{y}^{\beta}\partial_{t}^{m}F\left(  t,\cdot\right)  $ on the
left-hand side of the above display is square summable and the Parseval
relations (\ref{eq:P1})\ and (\ref{eq:P2}) hold. Further, $\partial_{y}%
^{\beta}\partial_{t}^{m}F=\partial^{\alpha}F\in L^{2}\left(  \mathbb{R}%
\times\mathbb{T}^{n}\right)  $ shows that formula (\ref{eq:P3}) is obtained by
integrating (\ref{eq:P1}) and invoking Fubini's theorem in the right-hand
side. Similarly, (\ref{eq:P4}) follows by integrating (\ref{eq:P2}), the use
of Fubini's theorem being allowed due to the estimate%
\begin{align*}
& \sum_{\xi\in\mathbb{Z}^{n}}\xi^{2\beta}\int\limits_{\mathbb{R}}\left\vert
\frac{d^{m}}{dt^{m}}\widehat{F}_{\xi}\left(  t\right)  \,\frac{d^{m}}{dt^{m}%
}\overline{\widehat{G}_{\xi}\left(  t\right)  }\right\vert dt\\
& \leq\left(  \sum_{\xi\in\mathbb{Z}^{n}}\xi^{2\beta}\int\limits_{\mathbb{R}%
}\left\vert \frac{d^{m}}{dt^{m}}\widehat{F}_{\xi}\left(  t\right)  \right\vert
^{2}dt\right)  ^{1/2}\left(  \sum_{\xi\in\mathbb{Z}^{n}}\xi^{2\beta}%
\int\limits_{\mathbb{R}}\left\vert \frac{d^{m}}{dt^{m}}\widehat{G}_{\xi
}\left(  t\right)  \right\vert ^{2}dt\right)  ^{1/2}\\
& =\left(  \int\limits_{\mathbb{R}}\int\limits_{\mathbb{T}^{n}}\left\vert
\partial^{\alpha}F\left(  t,y\right)  \right\vert ^{2}dydt\right)
^{1/2}\left(  \int\limits_{\mathbb{R}}\int\limits_{\mathbb{T}^{n}}\left\vert
\partial^{\alpha}G\left(  t,y\right)  \right\vert ^{2}dydt\right)
^{1/2}<\infty,
\end{align*}
which is a consequence of two Schwarz inequalities and formula (\ref{eq:P3}).
\endproof \medskip

\noindent\textbf{Proof of Theorem \ref{thm:poly-FI}.} Due to the absolute
convergence of the integrals and the series of formula (\ref{eq:P4}) with
$F:=S$, we have%
\begin{align}
\left\langle S,G\right\rangle _{\mathcal{B}_{p}}  & =\int\limits_{\mathbb{R}%
}\int\limits_{\mathbb{T}^{n}}\sum_{\left\vert \alpha\right\vert =p}\frac
{p!}{\alpha!}\partial^{\alpha}S\left(  t,y\right)  \,\overline{\partial
^{\alpha}G\left(  t,y\right)  }\,dydt\nonumber\\
& =\sum_{\xi\in\mathbb{Z}^{n}}\sum_{\left\vert \alpha\right\vert =p}\frac
{p!}{\alpha!}\xi^{2\beta}\int\limits_{\mathbb{R}}\frac{d^{m}}{dt^{m}}%
\widehat{S}_{\xi}\left(  t\right)  \,\frac{d^{m}}{dt^{m}}\overline{\widehat
{G}_{\xi}\left(  t\right)  }\,dt,\nonumber
\end{align}
where $\alpha=\left(  m,\beta\right)  $, as in Lemma~\ref{le:Parseval-Fubini}.
Next, replacing $\sum_{\left\vert \alpha\right\vert =p}$ by the double sum
$\sum_{m=0}^{p}\sum_{\left\vert \beta\right\vert =p-m}$ and using
$\alpha!=m!\beta!$ and the multinomial identity%
\[
\sum_{\left\vert \beta\right\vert =p-m}\frac{\left(  p-m\right)  !}{\beta!}%
\xi^{2\beta}=\left\vert \xi\right\vert ^{2\left(  p-m\right)  },
\]
we obtain%
\[
\left\langle S,G\right\rangle _{\mathcal{B}_{p}}=\sum_{\xi\in\mathbb{Z}^{n}%
}J_{\xi},
\]
where%
\[
J_{\xi}:=\sum_{m=0}^{p}\left(
\begin{array}
[c]{c}%
p\\
m
\end{array}
\right)  \left\vert \xi\right\vert ^{2\left(  p-m\right)  }\int_{-\infty
}^{\infty}\frac{d^{m}}{dt^{m}}\widehat{S}_{\xi}\left(  t\right)  \,\frac
{d^{m}}{dt^{m}}\overline{\widehat{G}_{\xi}\left(  t\right)  }\,dt.
\]
We now claim that the hypotheses of Lemma \ref{le:identity} and
Theorem~\ref{thm:FI-natLspline} are verified for $\sigma:=\widehat{S}_{\xi}$
and $\psi:=\widehat{G}_{\xi}$. Indeed, we have seen in the Introduction that
$S\in\mathcal{S}_{p}\left(  \tau,n\right)  $ implies $\widehat{S}_{\xi}%
\in\mathcal{S}_{p,\xi}\left(  \tau\right)  $, $\forall\xi\in\mathbb{Z}^{n}$.
On the other hand, $G\in\mathcal{B}_{p}$ implies $\widehat{G}_{\xi}\in
C^{p}\left(  \mathbb{R}\right)  $ and, by (\ref{eq:P3}), $\widehat{G}_{\xi
}^{\left(  m\right)  }\in L^{2}\left(  \mathbb{R}\right)  $, $\forall
m\in\left\{  0,\ldots,p\right\}  $, $\forall\xi\in\mathbb{Z}^{n}%
\backslash\left\{  0\right\}  $. The remaining condition $\widehat{G}_{\xi
}\left(  t_{j}\right)  =0$, $\forall j\in\left\{  0,\ldots,N\right\}  $
follows from $G\left(  t_{j},y\right)  =0$, $\forall j\in\left\{
0,\ldots,N\right\}  $, $\forall y\in\mathbb{T}^{n}$. Therefore Lemma
\ref{le:identity} and Theorem~\ref{thm:FI-natLspline} show that%
\[
J_{\xi}=\int_{-\infty}^{\infty}\left(  \frac{d}{dt}-\left\vert \xi\right\vert
\right)  ^{p}\widehat{S}_{\xi}\left(  t\right)  \left(  \frac{d}%
{dt}-\left\vert \xi\right\vert \right)  ^{p}\overline{\widehat{G}_{\xi}\left(
t\right)  }\,dt=0,\quad\forall\xi\in\mathbb{Z}^{n},
\]
which completes the proof. \endproof \medskip

\noindent\textbf{Proof of Theorem \ref{thm:var-p}.} This follows in a standard
way now by setting $G:=F-S$ in Theorem~\ref{thm:poly-FI}, to obtain
$\left\langle S,F-S\right\rangle _{\mathcal{B}_{p}}=0$, hence $0\leq\left\Vert
F-S\right\Vert _{\mathcal{B}_{p}}^{2}=\left\Vert F\right\Vert _{\mathcal{B}%
_{p}}^{2}-\left\Vert S\right\Vert _{\mathcal{B}_{p}}^{2}$. Equality is only
possible if $\left\Vert F-S\right\Vert _{\mathcal{B}_{p}}=0$, i.e.\ if $F-S$
is in the null space of the seminorm (\ref{eq:Duchon}). Since $p\geq2$, this
amounts to $\partial^{\alpha}\left(  F-S\right)  \equiv0$ for all
multi-indices $\alpha\in\mathbb{Z}_{+}^{n+1}$ with $\left\vert \alpha
\right\vert =p$. Taking $\alpha=\left(  0,\beta\right)  $ with $\beta
\in\mathbb{Z}_{+}^{n}$, it follows that the Fourier coefficients of $F-S$ with
respect to the last $n$ variables satisfy $\widehat{\left(  F-S\right)  }%
_{\xi}\left(  t\right)  =0$ for all $\xi\in\mathbb{Z}^{n}\backslash\left\{
0\right\}  $, $t\in\mathbb{R}$. Further, taking $\alpha=\left(  p,0\right)  $,
we find that the remaining Fourier coefficient for $\xi=0$ satisfies
$\frac{d^{p}}{dt^{p}}\widehat{\left(  F-S\right)  }_{0}\left(  t\right)  =0$,
$\forall t\in\mathbb{R}$. Hence $\left(  F-S\right)  \left(  t,y\right)  $ is
simply a polynomial of degree at most $p-1$ in $t$. Since $N\geq p-1$ and
$\left(  F-S\right)  \left(  t_{j},y\right)  =0$ for all $j\in\left\{
0,1,\ldots,N\right\}  $ and$\ y\in\mathbb{T}^{n}$, we obtain $F-S\equiv0$,
Q.E.D. \endproof \medskip

\noindent\textbf{Acknowledgements.} Thanks are due to Ognyan Kounchev
(Bulgarian Acade\-my of Sciences) and Michael Johnson (Kuwait University) for
stimulating discussions on the subject of this paper.

\end{document}